\documentclass[12pt,reqno]{amsart}

\def\b{\beta}
\def\g{\gamma}
\def\d{\delta}
\def\e{\epsilon}

\def\t{\theta}

\def\l{\lambda}
\def\m{\mu}

\def\x{\xi}

\def\o{\omega}

\def\D{\Delta}

\numberwithin{equation}{section}

\newtheorem {Theorem} 			{Theorem}
\newtheorem {varTheorem}                {Theorem}
\newenvironment {Theorem'}
        {\begin{varTheorem}{\hspace{-3.5mm}}{\bf '}{\hspace{3.5mm}}}
        {\end{varTheorem}}
\newtheorem {RefTheorem}[equation]     	{Theorem}  	
\newtheorem {Lemma}[equation]     	{Lemma}  	
\newtheorem {Proposition}[equation]	{Proposition}  

\newtheorem {Corollary}			{Corollary}

\theoremstyle{definition}
\newtheorem {Definition}[equation]{Definition}

\newtheorem {Example}[equation]		{Example}

\usepackage[cmtip,matrix,arrow]{xy}
\newdir{((}{{\lhook\kern-.1em}}
\newdir{))}{{\rhook\kern-.1em}}
\newdir{ >}{{}*!/-5pt/@{>}}

\def    \ker	{\operatorname{ker}}
\def    \im	{\operatorname{im}}

\def	\inv	{^{-1}}
\def	\to	{\longrightarrow}

\def	\t	{{\frak t}}
\def	\ft	{{\frak t}}

\def\lie {\rm Lie}

\def	\pr	{\operatorname{pr}}

\begin{document}
\title[On the cohomology rings of Hamiltonian T-spaces]
{On the cohomology rings of Hamiltonian T-spaces}

\author{Susan Tolman}
\author{Jonathan Weitsman}
\thanks{S. Tolman was partially supported by an NSF Mathematical Sciences
Postdoctoral Research Fellowship, by an NSF grant, and by an
Alfred P. Sloan Foundation Fellowship.
J. Weitsman was supported in part by NSF grant DMS 94/03567, by NSF 
Young Investigator grant DMS 94/57821, and
by an Alfred P. Sloan Foundation Fellowship.}

\address{Department of Mathematics,  University of Illinois at
Urbana-Champaign, 
Urbana, IL 61801}
\email{stolman@math.uiuc.edu}

\address{Department of Mathematics, University of California, Santa Cruz,
CA 95064}
\email{weitsman@cats.ucsc.edu}
\thanks{\today}

\begin{abstract}
Let $M$ be a symplectic manifold equipped with a Hamiltonian action of 
a  torus $T$.  Let $F$ denote the fixed point set of the $T$-action
and let $i:F\hookrightarrow M$ denote the inclusion.  By a theorem of
F. Kirwan \cite{K} the induced map $i^*:H_T^*(M)
\to H_T^*(F)$ in equivariant cohomology is an injection.  We give a simple
proof of a formula of Goresky-Kottwitz-MacPherson \cite{GKM}
for the image of the map $i^*$.
\end{abstract}

\maketitle



\section{Introduction}

The classification of manifolds equipped with group actions presents
difficulties beyond those inherent in the study of manifolds per se.
Even the basic questions---What is the equivariant cohomology ring of
the manifold? What can be said about the fixed manifolds of the group
action?  What is the cohomology ring of the quotient?---turn out to be
delicate and involved.

Much more is known in the case of a symplectic manifold $(M^{2m},\o)$ equipped
with a Hamiltonian action of a  torus $T$.  Let $H_T^*(M)$  denote
the rational equivariant cohomology ring of $M$.

The following theorem of F. Kirwan relates
the equivariant cohomology of $M$  with the equivariant
cohomology of its fixed point set:

\begin{RefTheorem}
\label{injectivity}(Kirwan)\cite{K}
Let a  torus $T$ 
act on a compact symplectic manifold $(M,\omega)$  
in a Hamiltonian fashion.
Denote the fixed point set of the  action by $F$.
The natural inclusion 
$i:F\to M$
of the fixed point set in the 
manifold   induces an  injection
$i^*:H^*_T(M)\hookrightarrow H^*_T(F)$.
\end{RefTheorem}

A related result, also due to Kirwan, relates the equivariant
cohomology of $M$ and the cohomology
of the  symplectic quotient of $M$.  

\begin{RefTheorem}
\label{surjectivity}
(Kirwan)\cite{K}
Let a  torus $T$ 
act on a compact symplectic manifold $(M,\omega)$  
with moment map $\mu: M \to \ft^*=\lie(T)^*$. 
Suppose that $\xi \in \ft^*$ is a regular value of $\mu$
and let $M_\xi := \mu\inv(\xi) /T$ be the
symplectic quotient of $M$.
The inclusion map $K:\mu\inv (\xi) \to M$ induces a surjection
$K^*: H^*_T(M) \to H^*_T(\mu\inv(\xi))=H^*(M_\xi)$.
\end{RefTheorem}

These two results give rise to two natural questions.  What is the image of
$i^*$?  And what is the kernel of $K^*$?  The purpose of this paper is to
address the first question.  In a companion paper \cite{TW},
we answer the second.

The description we give of the image of the injective map $i^*$ is due to
Goresky, Kottwitz, and MacPherson \cite{GKM}, using the work of Chang
and Skjelbred (\cite{CS}; see also \cite{AP,BV,H}).

\begin{Definition}
\label{1-skel}
Let $N \subset M$ denote the subset of $M$ consisting of points whose 
$T$-orbits are one-dimensional; i.e.,

$$N := \{p\in M \mid T \cdot p \sim S^1\}$$
\end{Definition}

By the local normal form theorem, each connected component $N_\alpha$ of $N$
is an open symplectic manifold whose closure $\overline{N_\alpha}$ is a
compact, symplectic submanifold of $M$. The restriction of $\mu$ to
$\overline{N_\alpha}$ is a moment map for the restricted torus action.
Furthermore the closure of $ N$ is given by 
$\overline{N} = N \cup F$;
this is referred to as the {\bf one-skeleton} of $M$.

\begin{Example}
Consider the case where $M^{2m}$ is a $2m$ dimensional toric
variety, equipped with the appropriate Hamiltonian action of a 
torus $T$ of rank $m$. 
The image of the moment map is the moment polytope $\D=\mu(M)$.
Let $v(\D)$ denote the union of of the vertices
of $\D$, and $ e(\D)$
the union of the interior of the edges.  Then $F= \mu\inv(v(\D))$, while
$N= \mu\inv(e(\D))$. 
\end{Example}

The main result of this paper is the following

\begin{Theorem}\label{main}
Let $(M,\omega)$ be a compact symplectic manifold equipped with a Hamiltonian
action of a  torus $T$.  
Let $F$ be the fixed point set and let $\overline{N}$ be the one-skeleton.
Let $i: F \to M$ and $ j:  F \to \overline{N} $  be the  natural
inclusions, and
$ H^*_T(M) \stackrel{i^*}{\to}  H^*_T(F) $  and
$ H^*_T(\overline{N}) \stackrel{j^*}{\to}  H^*_T(F) $ 
be the pull-back maps in equivariant cohomology.
Then the images of $i^*$ and $j^*$ are the same.
\end{Theorem}

This theorem is proved in considerable generality in \cite{GKM}.  The purpose
of this paper is to give a simple proof of this result in the symplectic
setting, which will enable us to obtain a description of the cohomology
ring of $M$ in a form that will make the structure of the map $K^*$ to the
cohomology ring of the symplectic quotient transparent.
Our methods
should yield similar statements in integral cohomology.  Additionally,
our methods give an algorithm for turning this description of the
cohomology ring into an explicit set of generators and relations. 

Before outlining the proof of Theorem \ref{main},  let us consider 
a few special examples.

\begin{Example}
If the torus $T$ is one-dimensional, $\overline{N}=M$, so 
the theorem is obviously true but trivial. 
\end{Example}

\begin{Example}
Suppose
that the closure $\overline{N_i}$ of each of the components of $N$ is a copy
of the two-sphere $P^1$.  For each such component there exists a
corank-one subgroup $K_i \subset T$ which acts trivially on $\overline{N_i}$.
The quotient $T/K_i$ is  isomorphic to $S^1$, and the corresponding action
on $P^1$ must be the usual action, so that $H^*_T(\overline{N_i})$ is 
given as follows.  Let $\g_i:K_i\to T$ denote the inclusion, and let 
$\g_i^*:H^*_T({\rm pt}) \to H^*_{K_i}({\rm pt})$ denote the 
induced map in equivariant cohomology.  For each $i$, the set $\overline{N_i}
\cap F$ consists of two points
$n_i, s_i$; and the image of $H^*_T(\overline{N_i})$ in $H^*_T(\overline{N_i}
\cap F)=H^*_T(n_i) \oplus H^*_T(s_i)$ consists of those elements 
$(a,b) \in H^*_T(n_i) \oplus H^*_T(s_i)$ such that $\g_i^*a=\g_i^*b$.
Let the fixed points of the $T$-action be given by $F_i$, $i=1,\dots,N$.  Then the
image of $H^*_T(M)$ in $H^*_T(F)= \bigoplus H^*_T(F_i)$ consists of
\begin{equation}\label{gkmcase}
(a_1,\dots, a_N) \in H^*_T(F_1)\oplus\dots H^*_T(F_N)
\end{equation}
such that, for each $\overline{N_i}$,
$$\gamma_i^*a_{n_i} = \gamma_i^*a_{s_i}$$
(See \cite{GKM}, \cite{GSbook}).
This gives a completely combinatorial algorithm for computing $H^*_T(M)$.
\end{Example}

The main tool needed to prove Theorem 1, as well as the injectivity and
surjectivity theorems \ref{injectivity}, \ref{surjectivity},
is a repeated use of equivariant Morse theory.
The key fact in all these cases is that components of the moment
map $\mu$ give equivariantly self-perfecting
Morse functions whose critical set is precisely $F$.
As the same is true for each of the $\overline{N_i}$'s, similar
statements can be made for the one-skeleton $\overline{N}.$

These self-perfecting Morse functions give us a very useful
way of constructing the equivariant cohomology ring of $M$ from
the cohomology rings of the fixed manifolds $F_i$:  roughly speaking,
the contribution
of each fixed manifold to the cohomology ring of $M$ consists
of those classes in $H^*_T(M)$ which vanish on all fixed points ``below''
$F_i$, and whose value on $F_i$ is a multiple of the downward Euler
class of the Morse flow.  As a similar statement can be made about
the cohomology ring of $\overline{N}$, we may compare the
images of $H^*_T(M)$ and $H^*_T(\overline{N})$ to prove our result.

\section{Morse Theory and the Moment Map}

In this section we state several results which will be the key steps
in the proof of Theorem \ref{main}.  Among them is Kirwan's injectivity
theorem (of which we supply a proof).  All of these results follow 
directly from the equivariant Morse complex associated to the choice
of a Hamiltonian as a Morse function.  We note that several of the
results of this section have analogs in integral cohomology; however 
we are only concerned with rational cohomology in this paper.

Let us recall our set-up.  Let $(M,\o)$ be a compact symplectic manifold with a
moment map $\m$ for the action of a  torus $T^n$.  Let $F\subset M$ 
denote the fixed point set of the torus action.  Given a generic element
$\x \in \frak{t}$, the function $f = \langle \m, \x \rangle : M \to {\bf R}$
is a Morse function on $M$ whose critical set coincides with the fixed
point set $F$.

Let us consider the fundamental exact sequence corresponding to the Morse
function $f$.  Let us denote by $C$ the critical set of $f$, and choose
$c\in C$.  We may assume that an interval $[c-\e, c+ \e]$ contains no critical
values of $f$ other than $c$.  Let $M^+_c = f^{-1} (-\infty, c + \e)$, 
$M^-_c= f^{-1} (-\infty, c-\e)$.  Then we have the following lemma, which 
is the main technical fact behind our results:

\begin{Proposition} \label{exact}
Let a  torus $T$ act on a compact  manifold $M$ with moment
map $\mu : M \to \t^*$.  Given $\xi \in \t$, choose a critical value
$c$ of the projection $f := \mu^\xi$. 
Let $F$ be the set of fixed points,  and let $F_c$ be the 
component of $F$ with value $c$.

The long exact sequence in equivariant cohomology for the pair
$(M^+_c, M^-_c)$ splits into short exact sequences:
\begin{equation} \label{exact:1}
0 \to H^*_T (M^+_c,M^-_c) {\to}H^*_T(M^+_c)
\stackrel{k^*}{\to} H^*_T(M^-_c) {\to} 0.\end{equation}
Moreover, the restriction from 
$H^*_T(M^+_c)$ to $H^*_T(F_c)$  induces an isomorphism
from the kernel of $k^*$ to those classes in $H^*_T(F_c)$
which are  multiples of
$e_c$, the equivariant Euler class  of the negative normal bundle
of $F_c$.
\end{Proposition}

Proof:  By our assumptions, $f$ is a Morse function, and there is a unique
critical value of $f$ contained in the interval $[c-\e,c+\e]$.  Denote the
corresponding critical manifold by $F_c$, and the negative disc and sphere
bundles of $F_c$ by $D_c$, $S_c$ respectively. 
The pair $(M_c^+, M_c^-)$ can be retracted onto the pair $(D_c, S_c)$, so there
is an isomorphism 

\begin{equation}\label{retr}
H^*_T(M_c^+,M_c^-)\cong H^*_T(D_c,S_c)\end{equation}
By the Thom isomorphism theorem, we have

\begin{equation}\label{thom}
H^*_T(D_c,S_c) \cong  H^{*-\l_c}_T(D_c) \  = H^{*-\l_c}_T(F_c)\end{equation}
where $\l_c$ is the Morse index of the critical manifold $F_c$; so we
obtain a commutative diagram

\begin{equation}\label{inj}
\vcenter{\xymatrix{
\ar[r] &
H^*_T(M_c^+,M_c^-) 
\ar[r]^{\g_c} \ar[d] & 
H^*_T(M_c^+) 
\ar[r]^{\b_c} \ar[d] & 
H^*_T(M_c^-)
\ar[r] &
\\
&
H^*_T(D_c,S_c)
\ar[r]^{\d_c} \ar[d]& 
H^*_T(D_c)
& 
\\
&
H^{*-\l_c}_T(D_c)
\ar[ur]_{\cup e_c} &
\\
}}
\end{equation}

\noindent where $e_c=e(D_c)$ is the equivariant Euler class of the bundle $D_c\to F_c$.
The cup product map $\cup e_c$ is {\em injective}; the same therefore
is true of the maps $\d_c$ and $\g_c$, proving the lemma.

The following corollary is then immediate:

\begin{Corollary}\label{perfection}(see \cite{AB1,AB2,K})
The function $f$ is an equivariantly perfect Morse function on $M$.
\end{Corollary}

Another application of this proposition is the proof of the following
theorem of Kirwan.

\begin{Theorem}\label{injectivity_noncomp}
Let a torus $T$ act on a  symplectic manifold $(M,\omega)$ with
proper bounded below  moment map $\mu: M \to \t^*$.  Let $i: F \to M $ denote
the natural inclusion of  the set $F$ of fixed points.
The pullback map $i^*:H^*_T(M) \to H^*_T(F)$ is injective.
\end{Theorem}

\begin{proof}
Order the critical values of $f$ as $c_1 < c_2 < \dots < c_n$.
The theorem obviously holds for $f^{-1}(-\infty, c_1) = \emptyset.$
Assume the proposition holds for the manifold $M^- := f^{-1}(-\infty, c_i)$.
We will show that it will hold for the manifold $M^+ 
:= f^{-1}(-\infty, c_{i+1})$;
the result follows then by induction.

By Lemma \ref{exact}, we have a map of short exact sequences


%
\begin{equation}\label{CD}
\vcenter{\xymatrix{
0
\ar[r] &
H^*_T(M^+,M^-)  
\ar[r] \ar[d]^{\simeq} & 
H^*_T(M^+) 
\ar[r] \ar[d]^{i_+^*} & 
H^*_T(M^-)
\ar[r]  \ar[d]^{i_-^*} &
0
\\
0
\ar[r] &
H^*_T(F_i)
\ar[r] &
H^*_T(F \cap M^+)
\ar[r] & 
H^*_T(F \cap M^-)
\ar[r] &
0,
\\
}}
\end{equation}
where $F_i$ denotes the critical set with value $c_i$.
By induction, the inclusion $i_-$ of $F \cap M^-$ into $M^-$ induces an
injection in cohomology.  By Proposition \ref{exact} the image of 
$H^*_T(M^+, M^-)$ in $H^*_T(M^+)$ is embedded injectively in $H^*_T(F_i)$.
The theorem then follows by diagram chasing.
\end{proof}

\section{Proof of the Main Theorem}

We are now ready to prove our main theorem: the 
restriction map to the fixed point set 
induces an isomorphism  from the equivariant cohomology of the original
Hamiltonian manifold to the image of the equivariant cohomology of the 
one-skeleton,  under its restriction map to the fixed point set. 
The key idea is to  use the tools developed in the last section
to  compare the graded rings associated to these images
using the filtration given by the Morse function
obtained by choosing a projection of the moment map.  The result will
then follow from the naturality of these objects, induction on
the critical points, and the injectivity of the restriction map.

Recall that the one-skeleton is given by
$$\overline{N}:=
\{p \in M \mid T \cdot p \mbox{ is one-dimensional or zero-dimensional}\}.$$
Clearly, the image of $i^*: H^*_T(M) \to H^*_T(F)$ is
a subset of the image of $j^*: H^*_T(\overline{N}) \to H^*_T(F)$. 
Therefore, $i^*$
induces a map, which we will also call $i^*$, from 
$H^*_T(M)$ to $\im j^* \subset H^*_T(F)$. 
By Theorem \ref{injectivity}, this map is injective.
Therefore, to prove the theorem, it suffices to show that this map
is surjective.  

On the level of the graded rings associated to the Morse filtration,
surjectivity will follow from comparing the proposition below with
Proposition~\ref{exact}.

\begin{Proposition}\label{onesk}
Let a  torus $T$ act on a compact symplectic manifold $M$ with
moment map $\mu : M \to \t^*$.
Given $\xi \in \t$, choose  a critical point $c$ of the
projection $f := \mu^\xi$. 
Let $F$ denote the fixed point set and let 
$F_c$ denote the component of $F$ with value $c$.
Define 
$F^- :=  F \cap f^{-1}(-\infty,c-\e)$ and
$\overline{N^+} :=  \overline{N} \cap  f^{-1}(-\infty,c+\e)$.

Let $\eta$ be a cohomology class in $H^*_T(\overline{N^+})$ 
which vanishes when restricted to $H^*_T(F^-)$.
Its restriction to $H^*_T(F_c)$ is a  multiple of $e_c=e(D_c)$, 
the equivariant Euler class of the downward normal bundle $D_c$
of $F_c$ (in $M$).
\end{Proposition}

\begin{proof}

Consider any component $N_\alpha$ of the set $N$
of one-dimensional orbits
such that the closure $\overline{N_\alpha}$ contains  $F_c$. 
The closure $\overline{N_\alpha}$ is a smooth $T$-invariant
symplectic manifold with moment map $\mu$.  
The class $\eta$ induces a cohomology class on 
$\overline{N_\alpha^+}  := \overline{N_\alpha} \cap f^{-1}(-\infty,
c+\epsilon)$ 
which vanishes when restricted to $\overline{N_\alpha} \cap F^-$,
and hence (by injectivity which we need to state in this version),
when restricted to 
$\overline{N_\alpha^-} := \overline{N_\alpha} \cap f^{-1}(-\infty, c-\epsilon)$.
Thus, by Proposition \ref{exact}, any element of the kernel of the natural map
$H^*_T(N_\alpha^+) \to H^*_T(N_\alpha^-)$ is, when restricted to
$H^*_T(F_c)$, 
a multiple of the equivariant Euler class $e_\alpha$
(here   $e_\alpha=e(D_c \cap \overline{N_\alpha})$
is the equivariant Euler class $e_c$ of the downward normal bundle $D_c \cap 
\overline{N_\alpha}$ of $F_c$
in $\overline{N_\alpha}$).

So the restriction of $\eta$ to $F_c$ is a multiple of the
equivariant Euler class of the downward normal bundle of $F_c$ in
$N_\alpha$.
Since this holds for each component $N_\alpha$,
and each of these  components  must have a different stabilizer, we
may apply Lemma \ref{lemma:4.2}  below. 
Therefore,
the class $\eta$ must be  multiple of the product of 
the equivariant  Euler classes the negative normal bundles to $F_c$
in all the components of $N$ whose closure contains $F_c$. 
But this is precisely the equivariant Euler class 
of the negative normal bundle to $F_c$ in $M$.
\end{proof}

\begin{Lemma}\label{lemma:4.2}
Let a torus $T$ act on a complex vector bundle $E$ over a manifold $F$,
so that the fixed set is precisely $F$.
Decompose $E$ into the direct sum of bundles $E_\alpha$,
where each $E_\alpha$ is acted on with a different weight $\alpha \in \t^*$.
Let $e_\alpha$ be the Euler class of the sub-bundle $E_\alpha$.

Then if $y \in H^*_T(F) $ is a  multiple of $e_\alpha$ for each $\alpha$,
then $y$ is a multiple of the product of the $e_\alpha$.
\end{Lemma}

\begin{proof}
Assume first that $F$ is a single point.

Let $\alpha \in \t^*$ be the weight with which $T$ acts on the
sub-bundle $E_\alpha$. Since $\alpha$ is a linear function
on $\t$, it lies naturally in $H^*(BT)$ = $H^*_T(F) = {\rm Sym(\t^*)}$,
the algebra of symmetric polynomials on $\t$.
The equivariant euler class of $E_\alpha$ is given by
$e_\alpha = \alpha^{n_\alpha}$, where $n_\alpha$ is the
complex dimension of $E_\alpha$.  The $\alpha$ are distinct
by assumption, and non-zero since no point not in the zero section is
fixed by $T$.  Therefore, the $e_\alpha$ are  pairwise relatively prime.
(Recall that every polynomial ring over  ${\bf Q}$ is a unique
factorization domain.)

More generally, since
$F$ is fixed by $T$, $H_T^{*}(F) = H^*(F) \otimes H^*(BT)$.  
Thus, $H_T^*(F)$ is bigraded.  In particular,  
given any integer $i$,
any cohomology class $a \in H_T^*(F)$ has a well-defined
component $a_i \in H^i(F) \otimes H^*(BT)$, and the sum
of all such components is $a$ itself; we will call $a_i$
{\bf the  component of $a$ with $F$-degree $i$}.

Note that the  component of $e_\alpha$ with $F$-degree $0$
is precisely $\alpha^{n_\alpha}$.  By the previous discussion,
these are non-zero and pairwise relatively prime.
Therefore it is enough to prove that
if $e$ and $f$ are two cohomology classes whose
components $e_0$ and $f_0$ with $F$-degree zero are
relatively prime, and if $e$ and $f$ both divide
$\alpha$, then so does $e \cdot f$.
We will prove this by induction.

We claim  that if 
$e ( f \cdot w + x) = f (e \cdot w  + y) $,
where the components of $x$ and $y$ with $F$-degree $i$
vanish for all $i < k$, then there exist $x', y'$ such that
$e ( f \cdot w + x') = f (e \cdot w + y')$,
and such that the components of $x'$ and $y'$ with $F$-degree $i$
vanish for all $i < k + 1$. 
To see this, compare the
component of $F$-degree $k$ on the two sides of the equation
$e ( f \cdot w + x) = f (e \cdot w + y)$.
Cancelling out terms which appear on
both sides, we get $e_0 x_k = f_0 y_k$.
Since $e_0$ and $f_0$ are relatively prime polynomials,
this shows that there exists $z_k$ such that
$x_k = f_0 z_k$ and $y_k = e_0 z_k$.
\end{proof}

We now proceed to prove  that the map $i^*: H^*_T(M) \to \im j^* \subset
H^*_T(F)$  is surjective.
We proceed, as usual, by induction:  

Consider  any  critical point $c$ of the projection $f := \mu^\xi$. 
Define $M^+ := f^{-1}(-\infty, c+\epsilon)$,
and $M^- := f^{-1}(-\infty, c-\epsilon)$  for any sufficiently
small $\e$.
Let $\overline{N^+} :=  \overline{N} \cap M^+$,
$\overline{N^-} :=  \overline{N} \cap  M^-$,
$F^+ :=  F \cap {M^+}$, and
$F^- :=  F \cap {M^-}$. 
Let $i^+: F^+ \to M^+$, $i^-: F^- \to M^-$, $j^+ :F^+ \to \overline{N^+}$ and
$j^- :F^- \to \overline{N^-}$ denote the corresponding inclusion maps. 
It is enough to assume that the induced map 
$i^{-*} : H^*_T(M^-) \to \im j^{-*} \subset H^*_T(F^-)$ 
is surjective, and to prove that the induced map
$i^{+*} : H^*_T(M^+) \to \im j^{+*} \subset H^*_T(F^+)$ 
is also surjective.

Since the images of $H^*_T(M^-)$ and $H^*_T(N^-)$ inside
$H^*_T(F^-)$ are the same, it follows that the natural
restriction map $r$ from $\im j^{+*} \subset H^*_T(F^+)$
to  $\im j^{-*}\subset H^*_T(F^-)$ is surjective.
Thus, taking the  exact sequence (\ref{exact:1}) of Proposition \ref{exact},
we have a map of short exact sequences

\begin{equation}\label{CD2}
\vcenter{\xymatrix{
0
\ar[r] &
H^*_T(M^+,M^-) 
\ar[r] \ar[d] & 
H^*_T(M^+) 
\ar[r] \ar[d]^{i_+^*} & 
H^*_T(M^-)
\ar[r]  \ar[d]^{i_-^*} &
0
\\
0
\ar[r] &
\ker r
\ar[r] &
\im j^{+*} \subset H^*_T(F^+)
\ar[r]^{r} & 
\im j^{-*} \subset H^*_T(F^-)
\ar[r] &
0,
\\
}}
\end{equation}

By our inductive assumption, $i^{-*}$ is surjective. 

By Proposition \ref{onesk},  every element in $\ker r$ is a multiple of
$e_c$, the equivariant Euler class of the negative normal bundle
of $F_c$, the  component of the fixed point set with value $c$. 
On the other hand, by Proposition \ref{exact}, every multiple of $e_c$
is in the image of the restriction $H^*_T(M^+,M^-)$ to $H^*_T(F_c)$
Thus,   the arrow from $H^*_T(M^+,M^-)$ to $\ker r$ is 
surjective too.  

The result follows by a diagram chase.

\section{Some comments about integral cohomology}

We close with some comments about integral cohomology.
Unfortunately,
the integer version of Theorem 1 is  not true in general.
In fact even the injectivity theorem \ref{injectivity} will
not hold in integral cohomology without some assumptions.
However, both injectivity and a version of Theorem 1 can be proved
over the integers where 
certain restrictions are placed on the allowable
stabilizer subgroups.  

Perhaps the  easiest example where Theorem 1 is
not true for integer cohomology is that of
$S^1 \times S^1$ acting on $S^2 \times S^2$,
with speed two on each sphere.\footnote
{If the reader is disturbed by the fact that this action is not
effective, she or he may tack on another couple of $S^2$'s spinning
at speed one.} 
Using the moment map for the diagonal action as our Morse function,
we see that every cohomology class which vanishes outside a neighborhood
of  the north pole $\times$ the north pole must be a
multiple of $4 x^2$ when restricted to that point.
In contrast, there exists a cohomology class on the
one-skeleton which vanishes outside this neighborhood
but is only a multiple of $2 x^2$.  Essentially, the problem
is that the weights are not relatively prime.  It is easy
to place a condition on the stabilizer groups at each fixed
point in a way that negates this possibility.  This is essentially
all that can go wrong, and a version of
Theorem 1 can be expected to hold if such an assumption of relative
primality is made.

\end{document}